\documentclass[11pt,a4paper]{article}
\usepackage{cite}
\usepackage{amsmath, amsthm, amssymb,slashed,upgreek}
\usepackage{ifpdf}
\ifpdf
  \usepackage[pdftex]{graphicx}
  \usepackage{epstopdf}
\else
  \usepackage[dvips]{graphicx}
\fi
\parskip=\baselineskip
\def\Bbb{\mathbb}
\def\Tr{{\rm Tr}}
\def\16{{\bf 16}}
\def\1{{\bf 1}}
\def\2{{\bf 2}}
\def\3{{\bf 3}}
\def\4{{\bf 4}}
\def\ad{{\mathrm{ad}}}
 
 \def\Spin{{\mathrm{Spin}}}

 \def\frak{\mathfrak}
\def\bar{\overline}

\def\R{{\Bbb{R}}}\def\Z{{\Bbb{Z}}}

\def\hat{\widehat}
\font\teneurm=eurm10 \font\seveneurm=eurm7 \font\fiveeurm=eurm5
\newfam\eurmfam
\textfont\eurmfam=\teneurm \scriptfont\eurmfam=\seveneurm
\scriptscriptfont\eurmfam=\fiveeurm

\font\teneusm=eusm10 \font\seveneusm=eusm7 \font\fiveeusm=eusm5
\newfam\eusmfam
\textfont\eusmfam=\teneusm \scriptfont\eusmfam=\seveneusm
\scriptscriptfont\eusmfam=\fiveeusm

\font\tencmmib=cmmib10 \skewchar\tencmmib='177
\font\sevencmmib=cmmib7 \skewchar\sevencmmib='177
\font\fivecmmib=cmmib5 \skewchar\fivecmmib='177
\newfam\cmmibfam
\textfont\cmmibfam=\tencmmib \scriptfont\cmmibfam=\sevencmmib
\scriptscriptfont\cmmibfam=\fivecmmib

\numberwithin{equation}{section}

\def\d{\mathrm d}
\def\C{{\Bbb C}}
\def\Z{{\Bbb Z}}
\def\A{{\mathcal A}}
\def\bar{\overline}
\def\W{{\mathcal W}}
\def\d{{\mathrm d}}
\def\O{{\mathcal O}}
\def\CP{\Bbb{CP}}
\def\R{{\Bbb R}}

\def\Z{{\Bbb Z}}

\def\L{{\mathcal L}}
\def\M{{\mathcal M}}
\def\t{{\mathfrak t}}
\def\A{{\mathcal A}}
\def\ad{{\mathrm{ad}}}
\def\d{{\mathrm d}}
\def\R{{\Bbb R}}
\def\Z{{\Bbb Z}}
\def\C{{\Bbb C}}
\def\Tr{{\mathrm{Tr}}}

\def\t{{\mathfrak t}}

\def\A{{\mathcal A}}
\def\V{{\mathcal V}}
\def\W{{\mathcal W}}
\def\O{{\mathcal O}}
\def\L{{\mathcal L}}
\def\D{{\mathcal D}}
\def\be{\begin{equation}}
\def\ee{\end{equation}}
\begin{document}
\begin{titlepage}
\begin{flushright}

\end{flushright}
\vskip 1.5in
\begin{center}
{\bf\Large{Two Lectures On  Gauge Theory}}
{\bf\Large{ and Khovanov Homology}}
\vskip
0.5cm {Edward Witten} \vskip 0.05in {\small{ \textit{School of
Natural Sciences, Institute for Advanced Study}\vskip -.4cm
{\textit{Einstein Drive, Princeton, NJ 08540 USA}}}
}
\end{center}
\vskip 0.5in
\baselineskip 16pt
\begin{abstract} In the first of these two lectures, I use a comparison to symplectic Khovanov homology to motivate the idea that the Jones polynomial
and Khovanov homology of knots
 can be defined by counting the solutions of certain elliptic partial differential equations in 4 or 5 dimensions.  The second lecture is devoted to a description
of the rather unusual boundary conditions by which these equations should be supplemented. An appendix describes some physical background. (Versions of these lectures
have been presented at various institutions including
the Simons Center at Stonybrook, the TSIMF conference center in Sanya, and also  Columbia University and the University of Pennsylvania.)  \end{abstract}
\date{March, 2016}
\end{titlepage}
\def\Hom{\mathrm{Hom}}

\def\U{{\mathcal U}}
\section{Lecture One}

The first physics-based proposal concerning Khovanov homology of knots was made by Gukov, Vafa, and Schwarz \cite{GVS},
who suggested that vector spaces associated to knots that had been introduced a few years earlier by Ooguri and Vafa \cite{OV} were related to
what appears in Khovanov homology.   
A number of years later, I re-expressed this type of construction in terms of gauge theory and the counting of solutions of PDE's
\cite{Witten}.  That is the story I will describe today.   Several previous lectures are
available  \cite{WittenOne,WittenTwo}  (the second of these may be a better starting point)
and I will take a different approach here.

\begin{figure}
 \begin{center}
   \includegraphics[width=3in]{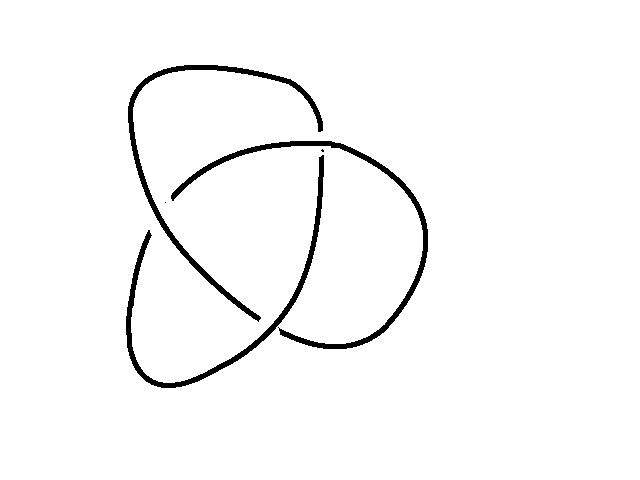}
 \end{center}
\caption{\small  A knot embedded in $\R^3$. \label{Knot}}
\end{figure}

In any event, the goal is to construct invariants of a knot embedded in $\R^3$ (fig. \ref{Knot}).
 In the simplest version, the invariants will be obtained by simply counting, with signs,
 the solutions of an equation.   The solutions will have an integer-valued\footnote{To be more precise, $P$ takes values in a $\Z$-torsor, rather than being canonically an integer.
 This is related to the framing anomaly of Chern-Simons theory. See Lecture 2.}  topological invariant $P$, 
 and if $a_n$ is the ``number'' (counted algebraically) of solutions with $P=n$,   then the Jones polynomial\footnote{In approaches based on quantum field theory, the natural normalization
 of the Jones polynomial of a knot or link in $\R^3$ is such that the Jones polynomial of the empty link is 1.  (The Jones polynomial is sometimes defined so that it equals 1 for an unknot
 rather than for the empty link.)   We normalize the argument $q$ of the Jones polynomial to be the instanton counting parameter, in a sense that will be explained later.  With
 this choice, the Jones polynomial of the unknot (with standard framing) is $q^{1/2}+q^{-1/2}$ and in general, for a knot with zero framing,  the exponents in eqn. (\ref{generating}) are half-integers.}
 of the knot will be
 \be\label{generating}J(q)=\sum_n a_n q^n. \ee
 
 To get Khovanov homology, this situation is supposed to be ``categorified,'' that is, we want for each $n$
 to define a complex of vector spaces whose Euler characteristic is $a_n$.     The only general situation that I know of in which one can naturally
 categorify the counting of solutions of an equation is the case that
 the equation whose solutions we are counting describes the critical points of some Morse function $h$.  We will be in this framework.
 Our
 equations will be partial differential equations or PDE's, so $h$  will be a Morse
 function on an infinite-dimensional space of functions, namely the functions that appear in the PDE.  The categorification will involve a middle-dimensional
 cohomology theory of the function space, analogous to Floer theory.  Let us put this aside
 for a moment and assume we are just trying to describe the  uncategorified theory, that is the Jones polynomial.
 
 The equations whose solutions I claim should be counted to define the Jones polynomial and ultimately Khovanov homology
 might look ad hoc if written down without an explanation of where they come from.    I could have started today's lecture by explaining
 the physical setup, but this might be unhelpful for some.   I decided instead to try a different approach of motivating the equations
by comparing to an established mathematical approach to Khovanov homology, namely symplectic Khovanov homology \cite{SS,Manolescu,AS}.
 
 Going all the way back to the original work of Vaughn Jones \cite{Jones}, most approaches to the Jones polynomial define an invariant
 in terms of some sort of presentation of a knot, for example a projection to a plane -- such as the projection used in drawing fig. \ref{Knot}.
     One defines something
 that is manifestly well-defined and explicitly computable once such a presentation is given.  What one defines is not obviously independent of the knot presentation,
 but turns out to be.   That step is where the magic is.  And there is always  some magic.
 
 An approach based on counting solutions of PDE's has the opposite advantages and drawbacks.  Topological invariance is potentially
 manifest (given certain generalities about elliptic PDE's and assuming compactness is under control), but it may not be clear how to calculate.   The ideal is to have manifest three- or (in the categorified case) four-dimensional symmetry
 together with a method of calculation.  How might this be achieved?
 
 I will suggest how to guess the right equations starting from a knowledge of symplectic Khovanov homology.  But in order to do this, we need
 to know something about a possible strategy to actually 
  count the solutions of an equation.  So I will begin by explaining what we would  do if we knew which equations we want
  to analyze, and this will help us in guessing the equations.
 
      \begin{figure}
 \begin{center}
   \includegraphics[width=3in]{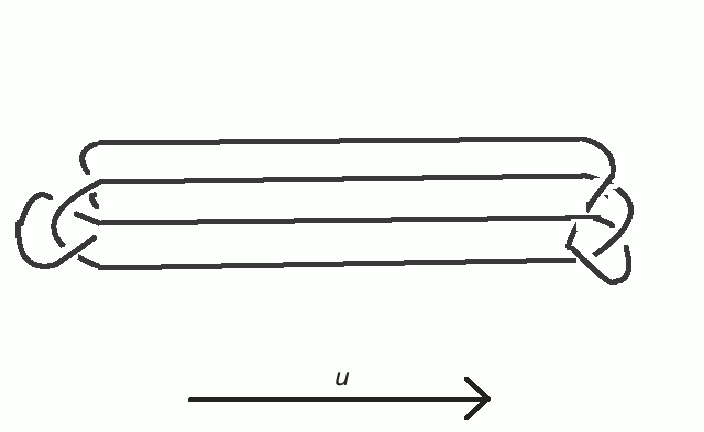}
 \end{center}
\caption{\small  A knot embedded in $\R^3$ and stretched in one direction. \label{longknot}}
\end{figure}
 There is a standard strategy, applicable to the present problem, for trying to count solutions of a PDE under suitable conditions.
 The original version was  the Atiyah-Floer conjecture concerning Floer homology of a three-manifold \cite{AtiyahFloer}.
Adapting their approach to the present problem, the idea is to stretch a knot in one direction, say the $u$ direction, as in fig. \ref{longknot}.
  Then one wants it to be the case that except near the ends, the solutions are independent of $u$. 
   This is not automatically the case and in \cite{GaiWit}, where this strategy was followed for the present problem, it was necessary to make a perturbation to a more generic system of equations to get to a situation in which this would be true.
  
\begin{figure}
 \begin{center}
   \includegraphics[width=3in]{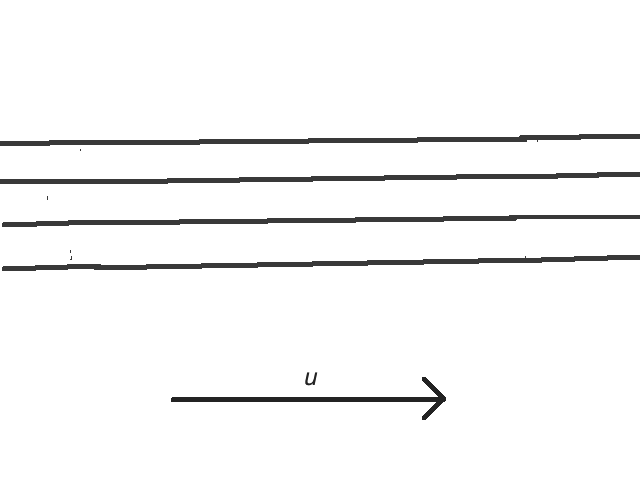}
 \end{center}
\caption{\small Infinite parallel strands parametrized by $u$, with $-\infty\leq u \leq \infty$. \label{longknotthree}}
\end{figure}

Given this, 
  we define a moduli space $\M$ of $u$-independent solutions.  We can think of these as the solutions in the presence of infinite parallel strands that run in the $u$ direction, as in fig. \ref{longknotthree}.

      \begin{figure}
 \begin{center}
   \includegraphics[width=3in]{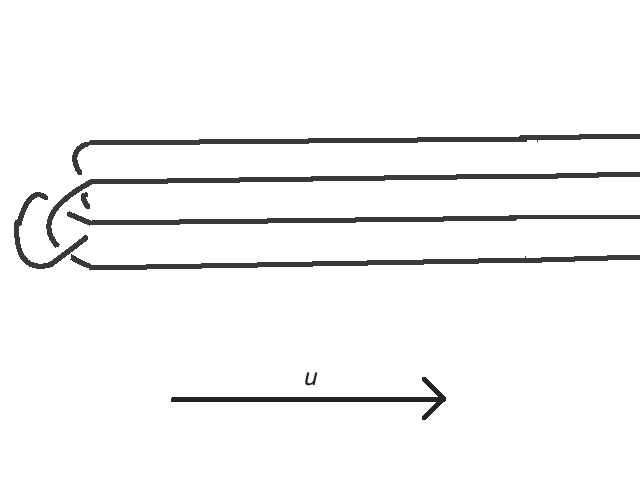}
 \end{center}
\caption{\small Semi-infinite strands that extend to $u=+\infty$. \label{longknottwo}}
\end{figure}
  
 Now as in fig. \ref{longknottwo} consider solutions in the presence of semi-infinite strands that extend to $u=+\infty$ or to $u=-\infty$ but not both.
 Let $\L_\ell$ and $\L_r$ be the moduli spaces of such solutions.  Thus a point in $\L_\ell$ represents a solution in a semi-infinite situation in which the strands terminate on the left (as drawn in fig. \ref{longknottwo}).
      Likewise $\L_r$ parametrizes solutions in the presence of semi-infinite strands that terminate on the right.
We assume that a solution in such a semi-infinite situation 
 is independent of $u$ for $u\to+\infty$ or $u\to-\infty$, respectively.  If this is so, then $\L_\ell$ and $\L_r$ come with natural maps to $\M$.  
 For simplicity in our terminology, we will assume that these maps are embeddings; this amounts to assuming that each solution in the interior in fig. \ref{longknot} can be extended over the left or over the right
  in at most one way.  This assumption is not necessary
 but makes the explanation simpler.

 The solutions for  a global knot like the one  in fig. \ref{longknot}
      can be understood as solutions in the middle that extend over both ends.  So the global solutions are intersection points
   of $\L_\ell$ and $\L_r$.      The integer $a_n$ that appears as a coefficient in the Jones polynomial is supposed to be the algebraic
   intersection number of $\L_\ell$ and $\L_r$:
   \be\label{intno} a_n=\L_\ell\cap \L_r.\ee
    (To be more exact, $a_n$ is this intersection number computed by counting only intersections with $P=n$.)

   In this language of intersections, categorification can happen if $\M$ is in a natural way a symplectic manifold and $\L_\ell$ and $\L_r$
   are Lagrangian submanifolds.    Then Floer cohomology -- i.e. the $A$-model or the Fukaya category -- of $\M$ gives a framework
   for categorification.    From the point of view of today's lecture, the reason that all this will happen is that, even before we stretched the knot
   to reduce to intersections in $\M$, the equations whose solutions we were counting are equations for critical points of some Morse function(al) $h$.
   
   In ``symplectic Khovanov homology,'' a version of such a story is developed for Khovanov homology (at least in a singly-graded version)
   with a very specific $\M$.    A description of this $\M$ that was proposed in \cite{K} (and exploited in a mirror version in \cite{CK}) 
   and which provided an important clue in my work is
   as follows.  $\M$ can be understood as a space of Hecke modifications.   Let me explain this concept.    Let $C$ be a Riemann surface
   and $E\to C$ a holomorphic $G_\C$ bundle over $C$, where $G_\C$ is some complex Lie group.   A Hecke modification of $E$ at a point $p\in C$
   is a holomorphic $G_\C$ bundle $E'\to C$ with an isomorphism to $E$ away from $p$:
   \be\label{maps}\varphi:E'|_{C\backslash p}\cong E|_{C\backslash p}. \ee
   
   For example, if $G_\C=\C^*$, the we can think of $E$ as a holomorphic line bundle $\L\to C$.    A holomorphic bundle $\L'$ that
   is isomorphic to $\L$ away from $p$ is
   \be\label{abn}\L'=\L(np)=\L\otimes \O(p)^n\ee
   for some integer $n$.   Here  $n$ can be thought of as a weight of the Langlands-GNO dual group of $\C^*$, which is another copy of
   $\C^*$.
   
   The reason that I write $G_\C$, making explicit that this is the complex form of the group, is that when we do gauge theory, the gauge group will be the
   compact real form and I will call this simply $G$.   In general, for any $G$, there is a corresponding Langlands-GNO dual group $G^\vee$,
   with complexification $G^\vee_\C$, such that Hecke modifications of a holomorphic  $G_\C$-bundle at a point $p\in C$ occur in families
   classified by dominant weights (or equivalently finite-dimensional representations) of $G^\vee_\C$ (or equivalently $G^\vee$).
   
   For example, if $G_\C=GL(2,\C)$, we can think of a $G_\C$-bundle $E\to C$ as a rank 2 complex vector bundle $E\to C$.
     The Langlanda-GNO dual group $G^\vee_\C$ is again $GL(2,\C)$, and a Hecke modification dual to the 2-dimensional representation of 
   $G^\vee_\C$ is as follows.   For some local decomposition $E\cong \O\oplus \O$ in a neighborhood of $p\in C$, one has
   $E'\cong \O(p)\oplus \O$.    The difference from the abelian case is that there is not just one Hecke modification of this type at $p$ but
   a whole {\it family} of them, arising from the choice of a subbundle $\O$ of $E$ that is going to be replaced by $\O(p)$.

     \begin{figure}
 \begin{center}
   \includegraphics[width=3in]{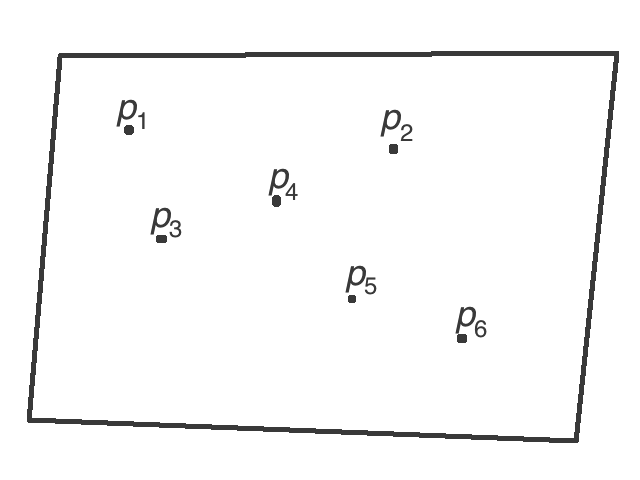}
 \end{center}
\caption{\small  A configuration of  points $p_i\in\R^2$, at which we are going to make Hecke modifications. \label{points}}
\end{figure}

 Because of this dependence, the Hecke modifications of this type at $p$ form a family, parametrized by $\CP^1$.  
Suppose we are given $2n$ points on $\C\cong \R^2$ at which we are going to make Hecke modifications of this type of a trivial bundle rank 2 complex
vector bundle $E\to \C$ (fig. \ref{points}).
The space of all such Hecke modifications would be a copy of $(\Bbb{CP}^1)^{2n}$, with one copy of $\Bbb{CP}^1$ at each point. 
However, there is a natural subvariety $\M\subset( \Bbb{CP}^1)^{2n}$ defined as follows.  One adds a point $\infty$ at infinity 
to compactify $\C$ to $\CP^1$, so
we are now making Hecke modifications of a trivial bundle $E=\O\oplus \O\to \Bbb{CP}^1$.   A point in $(\Bbb{CP}^1)^{2n }$ determines
a way to perform Hecke modifications at the points $p_1,p_2,\dots,p_{2n}$ to make a new bundle $E'$.   The space $\M$ is defined by requiring
that $E'\otimes \O(-n\infty)$ is trivial.   (If we were working in $PGL(2,\C)$ rather than $GL(2,\C)$, we would just say that $E'$ should
be trivial.) 

Symplectic Khovanov homology is constructed by considering intersections of Lagrangian submanifolds of the space $\M$ of multiple
Hecke modifications from a trivial bundle to itself. 
We want to reinterpret this in terms of gauge theory PDE's.

In my work with Kapustin on gauge theory and geometric Langlands \cite{KW}, an important fact was that $\M$ can be realized as a moduli space of
solutions of a certain system of PDE's.    However, although $\M$ is defined in terms of bundles on a 2-manifold $\R^2\cong\C$, the PDE's
are in 3 dimensions -- on $\R^3$.    As a result of this, everything in the rest of the lecture will be in a dimension one more than
one might expect.   To describe the Jones polynomial -- an invariant of knots in 3-space -- we will count solutions of certain PDE's in 4 dimensions,
 and the categorified version -- Khovanov homology -- will involve PDE's in 5 dimensions.  

The 3-dimensional PDE's that we need are known as the Bogomolny equations.  They are equations, on an oriented three-dimensional Riemannian manifold $W_3$,
for a pair $A,\upphi$, where $A$
is a connection on a $G$-bundle $E\to W_3$,   and $\upphi$ is a section of $\ad(E)\to W_3$
(i.e. an adjoint-valued 0-form).   If $F=\d A+A\wedge A$ is the curvature of $A$, then the Bogomolny equations are 
\be\label{z}F=\star \d_A\upphi.\ee
(Here $\star$ is the Hodge star and $\d_A$ is the gauge-covariant extension of the exterior derivative.) 

The Bogomolny equations have many remarkable properties and we will focus on just one aspect.  We consider the Bogomolny equations on $W_3=\R\times C$ with $C$ a Riemann surface.   Any connection $A$ on a $G$-bundle $E\to C$ determines a holomorphic structure on $E$ (or more
exactly on its complexification):  one simply writes $\d_A=\bar\partial_A+\partial_A$ and uses $\bar\partial_A$ to define the complex structure.  
(In complex dimension 1, there is no integrability condition that must be obeyed by a $\bar\partial$ operator.)   So for any $y\in \R$,
by restricting $E\to \R\times C$ to $E\to \{y\}\times C$, we get a holomorphic bundle $E_y\to C$.   However, if the Bogomolny equations are
satisfied, $E_y$ is canonically independent of $y$.    Indeed, a consequence of the Bogomolny equations is that $\bar\partial_A$ is independent of $y$
up to conjugation.    If we parametrize $\R$ by $y$, then the Bogomolny equations imply that
\be\label{zop}\left[\frac{D}{Dy}-\mathrm{i} \upphi,\bar\partial_A\right]=0.\ee
Thus $\bar\partial_A$ is independent of $y$, up to a natural conjugation.

The Bogomolny equations admit solutions with  singularities at isolated points.  
To understand the basic picture, we take the three-manifold to be simply $\R^3$, and the gauge group to be  $U(1)$.  One fixes an integer $n$ and one observes that the Bogomolny equation has an exact solution
for any $x_0\in\R^3$:
\be\label{zup}\upphi=\frac{n}{2|\vec x-\vec x_0|},~~ F=\star \d\upphi.\ee
  I have only defined $F$ and not the connection $A$ whose curvature is $F$ or the line bundle $\L$ on which $A$ is  connection. Such
an $\L$ and $A$ exist (and are essentially unique) if and only if $n\in \Z$.  

     \begin{figure}
 \begin{center}
   \includegraphics[width=3.5in]{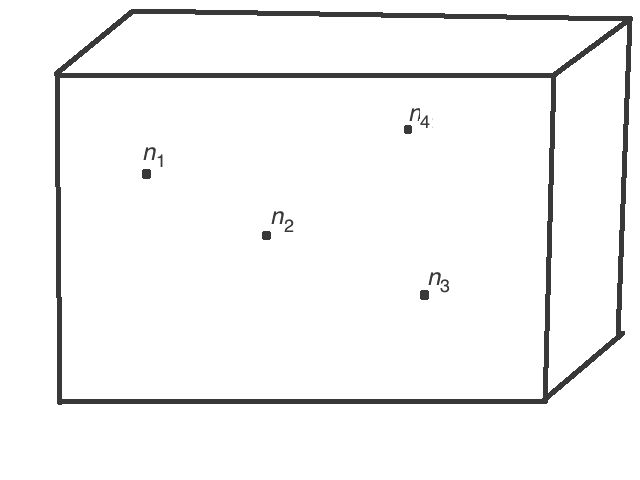}
 \end{center}
\caption{\small Points $y_i\times p_i\in\R^3$ labeled by weights $n_i$ of the group $U(1)$.  \label{bogsing}}
\end{figure}

For $G=U(1)$, since the Bogomolny equations are linear, they have a unique solution with singularities of this type labeled by specified integers $n_1,n_2,\dots$
at specified points $p_i\in \R^3$ (fig. \ref{bogsing}).  We simply take $\upphi=\sum_i \frac{n_i}{2|\vec x-\vec x_i|}$, $F=\star\d\upphi$.
   We assume that $\sum_in_i=0$, which ensures that $\upphi$ and the connection $A$ vanish at infinity faster than $1/|\vec x|$.
 
 Now pick a decomposition $\R^3=\R\times \R^2$, where we identify $\R^2$ as $\C$.   Suppose that the singularities are at $y_i\times p_i$,
 with $y_i\in\R$, $p_i\in\C$.
For each $y\notin\{y_1,\dots,y_n\}$, the indicated solution of the Bogomolny equations determines a holomorphic line bundle $\L_y\to \C$, and upon adding a point at infinity, this
naturally extends to $\L_y\to\CP^1$.  (Here we use the fact that $A$ vanishes at infinity faster than $1/|\vec x|$.)
$\L_y$ is independent of $y$  up to isomorphism as long as  $y$ is not equal to one of the $y_i$, but even when $y$ crosses one of the $y_i$, $\L_y$ is constant
when restricted to $\CP^1\backslash p_i$.    In crossing $y=y_i$, $\L_y$ undergoes a Hecke modification
\be\label{mapt}\L_y\to \L_y\otimes \O(p_i)^{n_i}.\ee
$\L_y$ is trivial for $y\to -\infty$ and for $y\to +\infty$ (again because the solution vanishes at infinity faster than $1/|\vec x|$).   The solution thus describes a sequence of Hecke modifications mapping the trivial bundle
to itself.

We can do something similar for any simple Lie group $G$.  (The underlying idea was introduced by 't Hooft in the late 1970's \cite{thooft} and is important in
physical applications of quantum gauge theory.)   Let $T$ be the maximal torus of $G$ and let $\frak t$ be its Lie algebra.  Pick a homomorphism
$\rho:\frak{u}(1)\to \frak t.$   Up to a Weyl transformation, such a $\rho$ is equivalent to a dominant weight of the dual group $G^\vee$, so
it corresponds to a representation $R^\vee$ of $G^\vee$.    We turn the singular solution (\ref{zup})  of the $U(1)$ Bogomolny equations that we already
used (more exactly, the special case of this solution with $n=1$) into a singular
solution for $G$ simply by
\be\label{realt}(A,\upphi)\to (\rho(A),\rho(\upphi)).\ee
  Then we look for solutions of the Bogomolny equations for $G$ with singularities of this type at specified points $y_i\times p_i\in \R^3$. 

     \begin{figure}
 \begin{center}
   \includegraphics[width=3.5in]{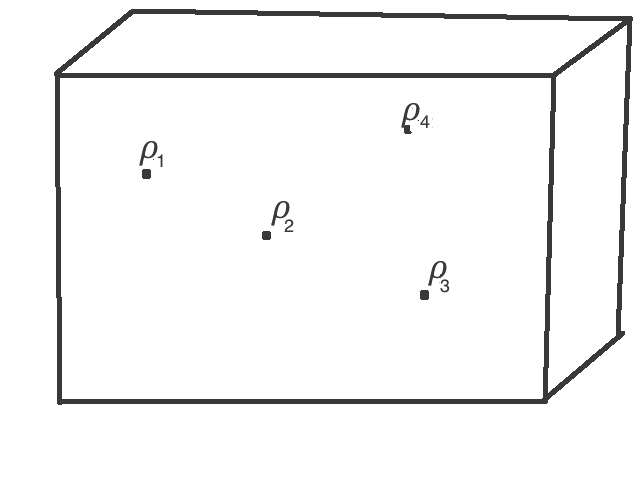}
 \end{center}
\caption{\small Points in $\R^3$ labeled by homomorphisms $\rho_i:\frak{u}(1)\to \frak t$, or equivalently by representations of the dual group.  \label{bogsing2}}
\end{figure}
 
The picture is the same as before  
  except that now (fig. \ref{bogsing2}) the points $y_i\times p_i$ are labeled by homomorphisms
  $\rho_i:\frak{u}(1)\to \frak t$, or in other words by representations $R_i^\vee$ of the dual group $G^\vee$, rather than by integers $n_i$. 
  Also, now we must specify that the solution should go to 0 at infinity faster than $1/r$ (for $U(1)$, this was automatic once we set $\sum_in_i=0$).  Given this, such a solution describes a sequence
  of Hecke modifications at $p_i$ of type $\rho_i$, mapping a trivial $G$-bundle $E\to \CP^1$ to itself. 
  
 The moduli space $\M$ of solutions of the Bogomolny equations on $\R^3$ with the indicated singularities and
vanishing at infinity faster than $1/r$ is actually a hyper-Kahler manifold, essentially first studied by P. Kronheimer in the 1980's. 
If we pick a decomposition $\R^3= \R\times \R^2$, this picks one of the complex structures on the hyper-Kahler manifold and in that
complex structure, $\M$ is the moduli space
 $\M_{p_1,\rho_1;p_2,\rho_2;\dots}$ of all Hecke modifications  of the indicated types at the indicated points, mapping a trivial bundle over $\CP^1$
 to itself.   
 
 This construction can be used to account for a number of properties of spaces of Hecke modifications, but for today we want to focus on the
  application to knot theory.  The reduction to $\M$ is supposed to result from stretching a knot in one direction, so 
  we want $\M$ to be the space of $u$-independent solutions of some equations, as suggested in fig. \ref{longknotthree}.
  We already described $\M$ via solutions of some PDE's on $\R^3$, so now we have to 
     think of $\M$ as a space of $u$-independent solutions on $\R^4=\R^3\times \R$, where the second factor is parametrized by $u$.

There actually are natural PDE's in four dimensions that work. They play a role in the gauge theory approach to
 geometric Langlands \cite{KW}, and are  sometimes called the KW equations.  They are equations for a pair $A,\phi$ where $A$ is a connection on $E\to Y_4$, $Y_4$ a four-manifold, and $\phi$ is a 1-form
on $Y_4$ valued in $\ad(E)$:
\be\label{meqn}F-\phi\wedge \phi=\star\d_A\phi,~~~\d_A\star\phi=0.\ee
In a special case $Y_4=W_3\times \R$, with $A$ a pullback from $W_3$ and $\phi =\upphi\, \d u$ (where $\upphi$ is a section of $\mathrm{ad}(E)$ and 
$u$ parametrizes the second factor in $Y_4$) these
equations reduce to the Bogomolny equations on $W_3$:
\be\label{redb}F=\star\d_A\upphi.\ee

\begin{figure}
 \begin{center}
   \includegraphics[width=3in]{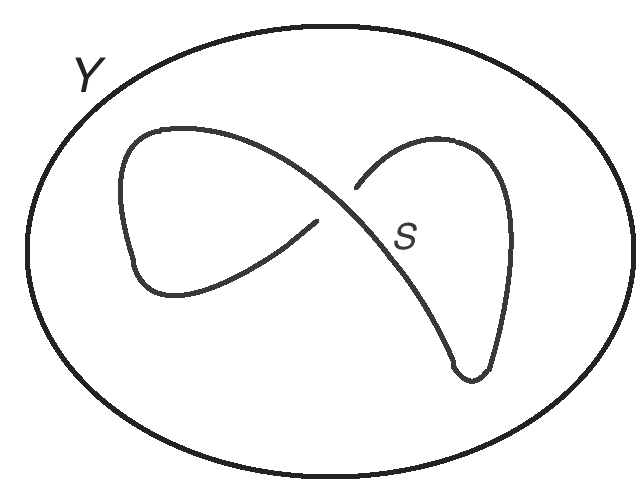}
 \end{center}
\caption{\small A four-manifold $Y$ with an embedded 1-manifold $S$ along which one specifies a desired singularity.  \label{LoopSing}}
\end{figure}

Therefore, the singular solution (\ref{realt}) of the Bogomolny equations that we have already studied can be lifted to a singular solution of the KW
equations, but now the singularity is along a line rather than a point. Of course, the singularity is still in codimension three. We view this solution as a model that tells us what sort of codimension three singularity to look
for in a more general situation.    If $Y_4$ is a 4-manifold and $S\subset Y_4$ is an embedded 1-manifold, labeled
by a homomorphism $\rho:\frak{u}(1)\to \frak{t}$ (or by a representation of $G^\vee$), then one can look for solutions of the KW equations with
a singularity  along $S$ associated to the given choice of $\rho$ (fig. \ref{LoopSing}).

If we specialize to the case that $Y_4=\R^3\times \R$, with   $S=\cup_i S_i$, and $S_i=p_i\times \R\subset \R^3\times \R$ ($p_i$ are
points in $\R^3$ and $\R$ is parametrized by $u$) 
     then the $u$-independent solutions of the KW equations are just the solutions of the Bogomolny equations on $\R^3$, with
     the chosen singularities.  So these solutions are parametrized by $\M$; and indeed one can show that these are all solutions of the KW
     equations in this situation with reasonable behavior at infinity.

So we have an elliptic PDE in four dimensions and we can specify in an interesting way what sort of singularity it should have on an embedded
circle $S\subset Y_4$.    But this sounds like a ridiculous framework for knot theory, because there is no knottedness of a 1-manifold in 
a 4-manifold!

To resolve this point, we have to explain what is involved in categorification.
Let us practice with an ordinary equation rather than a partial differential equation.   Suppose that we are on a finite-dimensional 
compact oriented manifold
$N$ with a real vector bundle $V\to N$ with rank$(V)$=dimension$(N)$.    Suppose also we are given a section $s$ of $V$. 
We can define an integer by counting, with multiplicities (and in particular with signs) the zeroes of $s$.    This integer is the Euler class
$\int_M\chi(V)$. 

In general as far as I know, there is no way to categorify the Euler class of a vector bundle.    However, suppose that $V=T^*N$
and that $s=\d h$ where $h$ is a Morse function.    Then the zeroes of $s$, which are critical points of $h$,
 have a natural ``categorification'' described in Morse homology.
One defines a complex $\mathcal V$ with a basis vector $\psi_p$ for each critical point $p$ of $h$.   The complex is $\Z$-graded by assigning
to $\psi_p$ the ``index'' of the critical point $p$, and it has a natural differential that is defined by counting gradient flow lines between different critical points.

    \begin{figure}
 \begin{center}
   \includegraphics[width=3in]{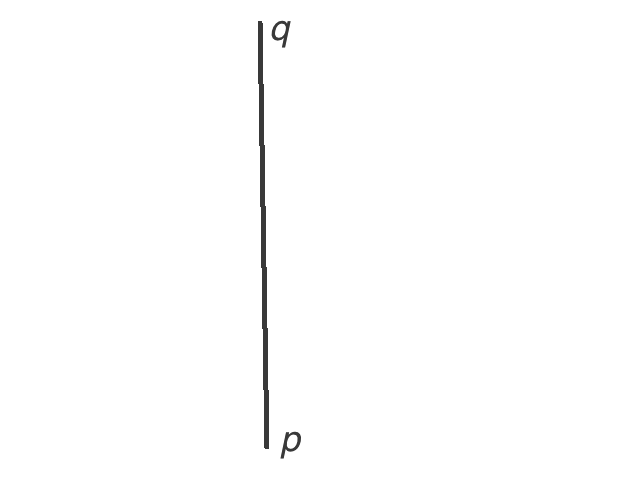}
 \end{center}
\caption{\small A flow from one critical point to another.  \label{flowg}}
\end{figure}
Concretely the differential is defined by
\be\label{diffdef}\d\psi_p=\sum_q n_{pq}\psi_q\ee
where the sum runs over all critical points $q$ whose Morse index exceeds by 1 that of p, and the integer $n_{pq}$ is defined by counting flows from $p$ to $q$ (fig. \ref{flowg}).
    A ``flow'' is a solution of the gradient flow equation
    \be\label{flowz}\frac{\d \vec x}{\d t}=-\vec\nabla h.\ee
    (To define this equation, one has to pick a Riemannian metric on the manifold $N$.  The complex that one gets is independent of the metric up to
    quasi-isomorphism.  One considers flows that start at $p$ at $t=-\infty$ and end at $q$ at $t=+\infty$.  Such flows come in one-parameter families related by time translations
    and $n_{pq}$ is the number of such families, counted algebraically.)
    
This tells us what we need in order to be able to categorify a problem of counting solutions of the KW equations.
We have to be able to write those equations  as equations for a critical point of a functional $\Gamma(A,\phi)$:
\be\label{mext}\frac{\delta \Gamma}{\delta A}=\frac{\delta \Gamma}{\delta \phi}=0.\ee
And the associated gradient flow equation, which will be a PDE in 5 dimensions on $X_5=\R\times Y_4$
\be\label{wext}\frac{\d A}{\d t}=-   \frac{\delta \Gamma}{\delta A},~~~~~\frac{\d \phi}{\d t}=-   \frac{\delta \Gamma}{\delta \phi}\ee
has to be elliptic, so that it will makes sense to try to count its solutions. 

Generically, it is not true that the KW equations on a manifold $Y_4$ are equations for a critical point of some functional.  
However, this is true if $Y_4=W_3\times \R$ for some $W_3$.   If singularities are present on an embedded 1-manifold $S\subset Y_4$
then there is a further condition:  The KW equations in this situation are equations for critical points of a functional if and only if $S$ is contained
in  a 3-manifold $W_3\times p$, with $p$ a point in $\R$. (For an explanation of ``why'' this  is true, see the appendix.)    So to make categorification possible, we have to be in the situation that leads to
knot theory:  $S$ is an embedded 1-manifold in a 3-manifold $W_3$.  Once this restriction is made, the five-dimensional flow equations exist and
are indeed elliptic.  (They were introduced
independently in \cite{Witten} and \cite{Haydys} and are sometimes called the HW equations.)

Naively, this leads to ``categorified'' knot invariants for any three-manifold $W_3$, but to justify this claim one needs some compactness
properties for solutions of the equations under consideration.    I suspect that a proper proof of these compactness properties may require
that the Ricci tensor of $W_3$ is nonnegative, a very restrictive condition. 

What I have described so far is supposed to correspond (for $W_3=\R^3$, $G=SO(3)$ and $\rho$ corresponding to the 2-dimensional representation
of $G^\vee=SU(2)$) to ``singly-graded Khovanov homology.''   It is singly-graded because the only grading I have mentioned
is the grading that is associated to the Morse index, or in other words to categorification.   In the mathematical theory, one says that singly-graded
Khovanov homology becomes trivial (it does not distinguish knots) if one ``decategorifies'' it and forgets the grading.    In the approach I have described,
this is true because in the uncategorified version, the embedded 1-manifold $S$ is just a 1-manifold in a 4-manifold $Y_4$ (it has no reason to
be embedded in the 3-manifold $W_3\times p$) so there is no knottedness.

The physical picture makes clear where the additional ``$q$''-grading of Khovanov homology would come from.   It is supposed to come
from the second Chern-class, integrated over the 4-manifold $Y_4$.  The second Chern class is the  invariant that I called $P$ at
the beginning of the lecture.  But for topological reasons, the second Chern class cannot be defined in the presence of  a  codimension three singularity of the type I have
described on an embedded one-manifold $S\subset Y_4$.
(Because of the singularity, $Y_4$ behaves as a noncompact four-manifold on which there is no topological invariant corresponding to the second Chern class
of a $G$ bundle.)
And therefore the construction as I have presented it so far has no $q$-grading.

     \begin{figure}
 \begin{center}
   \includegraphics[width=3in]{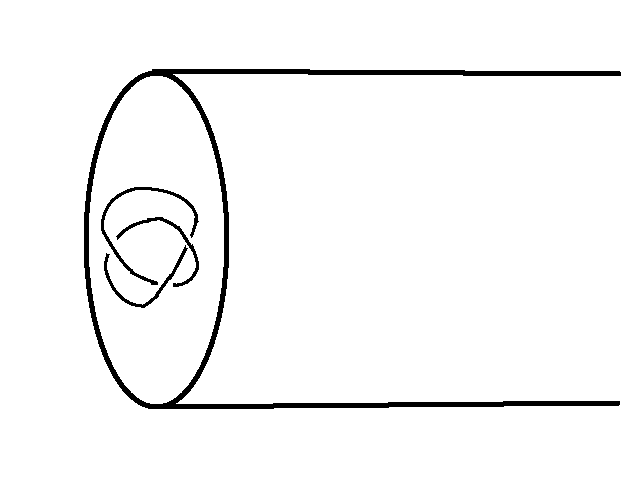}
 \end{center}
\caption{\small A four-manifold $Y_4$ with boundary, with a knot placed in its boundary.  \label{thimble}}
\end{figure}   
    The physical picture tells us what we have to do to get the $q$-grading:   $Y_4$ should be a manifold with boundary, with the knot
    placed in its boundary (fig. \ref{thimble}).
     The appropriate boundary condition will be the subject of Lecture Two and is such that the second Chern class
     can be defined.

 In \cite{GaiWit}, Gaiotto and I analyzed this problem  (in the uncategorified situation, meaning that we counted solutions in 4 dimensions, not 5, and for the simplest
 case of $G=SO(3)$) with the aim
 of showing directly, without referring to the physical picture, that the Jones polynomial is
 \be\label{pext}J(q)=\sum_na_nq^n\ee
 where $a_n$ is the number of solutions with second Chern class $n$.    The starting point was to stretch the knot in one direction,
 reducing to equations in one dimension less, as in fig. \ref{longknot}. It turns out that 
  the solutions in one dimension less that satisfy the boundary condition are related to a lot of interesting mathematical physics
 involving integrable systems, conformal field theory, geometric Langlands, and more.
What emerges is the ``vertex model'' construction of the Jones polynomial; the way it emerges is somewhat along the lines of work by Bigelow \cite{Bigelow} and  Lawrence \cite{Lawrence}.   
What our work added was a derivation of the vertex model from
 a starting point with manifest 3-dimensional symmetry.  The analog of this for the categorified theory   is expected to involve, in one version,
  a Fukaya-Seidel category with a certain superpotential.    The relevant model -- as well as a plausible variant that does not work -- has been explored and to a considerable
  extent understood in \cite{GalMoore}.

  \section{Lecture Two} 

In Lecture One, I explained that to define a $q$-grading in Khovanov homology, we have to be able to make sense of the second
Chern class of a solution of the KW equations on a four-manifold $Y_4$, in the presence of a knot.  As already explained, if the knot
is represented by a codimension 1 embedded submanifold $S\subset Y_4$, this will not work, because the singularity that we want to postulate
along $S$ does not allow the definition of a second Chern class as a topological invariant.  Instead, we embed the knot in the boundary of $Y_4$, as
 in fig. \ref{thimble}.
     The boundary condition that we use is subtle to describe, but has the property that the bundle is fixed on the boundary, so the second Chern class
     can be defined.

We could actually get the $q$-grading for any $Y_4$ with boundary, but to also allow categorification, we want more specifically
$Y_4=W\times \R_+$, where $\R_+$ is a half-line, parametrized by $y$. 
For the Jones polynomial and Khovanov homology, we further take $W=\R^3$. (More general choices of $W$ are certainly also interesting, but not much is known
about what to expect. See \cite{GPV}.)  So as sketched in fig. \ref{thimble},
$Y_4$ is $\R^3\times \R_+$  with the knot embedded in the boundary.

For $y\to \infty$, we ask for $A,\phi\to 0$.    For $y\to 0$, there is a subtle boundary condition which is one of the main points
of the theory.   Describing it is actually my main goal for today.  This boundary condition depends on the knot  $K$, and on the labeling
of $K$ (or of each component of a link $L$) by a representation $R^\vee$  of the Langlands or GNO dual group $G^\vee$.   That is the
only place that $K$ enters the setup.    

The desired boundary condition is an elliptic boundary condition though possibly an unexpected one.   Here, ``elliptic'' means that although
the definition may look unexpected, the resulting properties are similar to what one would get with more familiar elliptic boundary conditions
such as Dirichlet or Neumann.   For example, on a compact four-manifold, the linearization of the KW equation becomes a Fredholm operator and  has discrete spectrum.

With this boundary condition, which I will describe in some detail,
the restriction to the boundary $W\times \{y=0\}$ of the bundle $E$ and connection $A$ are specified.   As a result, one can define
a second Chern class
\be\label{zort}n=\frac{1}{8\pi^2} \int_{W\times \R_+}\Tr\, F\wedge F.\ee
However, because the bundle is fixed on the boundary but is not trivialized, this invariant really takes values in a $\Bbb Z$-torsor associated to framings of $W$ and $K$.    ``Torsor'' means that it is not
true that $n$ is an integer in a canonical way; rather the value of $n$ mod $\Z$ depends only on the boundary conditions and not on the specific gauge field
that satisfies them.   One can ``trivialize the torsor'' and make $n$ an integer by picking framings of $W$ and $K$.

To define a knot polynomial, one counts (with signs, in a standard way) the number $	a_n$ of solutions of the KW equations
with second Chern class $n$, and then one defines
\be\label{peft}J(q;K,R^\vee)=\sum_{n\in\Bbb Z}a_nq^n.\ee
  Compactness (not yet proved with the appropriate boundary conditions) of the solutions of the KW equations will mean
that there are only finitely many terms in the sum so that this is a Laurent polynomial.

As I explained in Lecture One, with this definition of the Jones polynomial, the ``categorification'' that leads to Khovanov homology is straightforward
in principle.    It arises because the KW equations can be ``lifted,'' in a certain sense, to certain elliptic differential equations in five
dimensions, and these equations can be interpreted as gradient flow equations.   But rather than say more about that today,
what I want to do is to describe the boundary condition that is needed for the four (or five) dimensional equations.    This boundary
condition is of a possibly somewhat unfamiliar type, and understanding it is essential for making progress with this subject.

The boundary
conditions have been studied in \cite{MW},  and were shown to be elliptic in the absence of a knot.    A paper is in progress on the case with a knot \cite{MWtwo},
and I will tell a little
about that case later.    I will carry out this discussion in the language of the four-dimensional equations, since going to five dimensions does not change much,
as was shown in \cite{MW}.

The boundary condition of interest is not a simple Dirichlet or Neumann boundary condition -- it is not defined by saying what fields or
derivatives of fields vanish along the boundary.    Rather, the boundary condition is defined by specifying a model solution of the KW
equations that has a singularity along the boundary, and saying that one only wants to consider solutions of the KW equation that
are asymptotic to this singular solution along the boundary.

The model solution is a solution on $\R^3\times \R_+$, where I will parametrize $\R^3$ by $x_1,x_2,x_3$ and $\R_+$ by
$y$.  There is a simple exact solution with $A=0$ and
\be\label{yef}\phi=\sum_{i=1}^3\frac{\t_i \cdot \d x_i}{y}, \ee
where $\t_i$ are elements of the Lie algebra $\frak g$ of $G$ that obey the $\frak{su}(2)$ commutation relations
\be\label{yyef}[\t_1,\t_2]=\t_3,~~\mathrm{and~cyclic~permutations}.\ee
  Thus the $\t_i$ are images of a standard basis of $\frak{su}(2)$ under some homomorphism $\rho:\frak{su}(2)\to \frak g$.
  Every $\rho$ leads to an interesting theory, but to get the Jones polynomial and Khovanov homology, we take $\rho$ to be
a principal embedding in the sense of Kostant.  (For $G=SU(2)$, this simply means that $\rho$ is the identity map
$\frak{su}(2)\to \frak g$. For $G=SU(N)$, it means that the $N$-dimensional representation of $\frak g$
is irreducible with respect to $\rho(\frak{su}(2))$.)

The solution I have just described is what I call the Nahm pole solution, since the relevant singularity was introduced long ago
by Nahm  in his work on magnetic monopoles \cite{N}.   That was in the context of ``Nahm's equation,'' which is an ordinary differential equation for three
$\frak g$-valued functions $\phi_1,\phi_2,\phi_3$ of a real variable $y$:
\be\label{zefa}\frac{\d\phi_1}{\d y}+[\phi_2,\phi_3]=0,~~\mathrm{and~cyclic~permutations}.\ee
The KW equations reduce to Nahm's equation if we drop the dependence on $\vec x$
and set $A=\phi_y=0$.
  On $\R^3\times\R_+$, the Nahm pole boundary condition
just says that a solution is supposed to be asymptotic to the Nahm pole solution for $y\to 0$.

To state the Nahm pole boundary condition on $M_4=W_3\times \R_+$, for a more general 3-manifold $W_3$, one needs to specify some
terms of $\O(1)$ in the solution (for $y\to 0$) as well as the singular terms of order $1/y$.   For $G=SU(2)$, one takes the $G$-bundle
$E$ on which we are solving the equations to be, when restricted to $W_3\times \{y=0\}$, the frame bundle of $W_3$, so that $\ad(E)=TW_3$. 
Then one takes $A$, restricted to the boundary, to be the Levi-Civita connection on $TW_3$. 
 With this choice of $E$, the formula  
\be\label{torf}\phi=\sum_{i=1}^3\frac{\t_i \cdot \d x_i}{y} \ee
makes sense (one can think of the numerator $\sum_i\t_i\cdot \d  x_i$ as stating the identification $\ad(E)\cong TW_3$).   One can show that this choice
of $(A,\phi)$ obeys the Nahm pole boundary condition up to an error of $\O(y)$, and the Nahm pole boundary condition simply says that
the solution should agree with what I have described up to $\O(y)$.  (One can  generalize this to the case that the metric of $M_4$
is not a product near the boundary.)

Showing that the Nahm pole boundary condition is elliptic is mostly an exercise in ``uniformly degenerate elliptic operators,'' but one
needs to know some specific facts about the KW equations.   The main thing that one needs to know is that if $\L$ is the linearization
of the KW equations on a half-space $\R^4_+$ around the Nahm pole solution, then $\L$ as an operator between appropriate Hilbert spaces of functions on $\R^4_+$
has no kernel or cokernel.     Actually one can show in an elementary way that $\L^\dagger =-N\L N^{-1}$ with an explicit matrix
$N$, so it suffices to show that there is no kernel. 

Much the same argument that proves this actually proves the following statement:    The only solution of the KW equations on
$\R^4_+$, approaching the Nahm pole solution for  $y\to 0$ and also for $\sqrt{\vec x^2+y^2}\to\infty$, is the Nahm pole solution, and moreover this solution
is ``transverse'' (ln expanding around it, the operator $\L$ has zero kernel and cokernel).  In terms of
Khovanov homology, this means that the Khovanov homology of the empty link is of rank 1.

Before trying to prove these vanishing results, I will explain a simpler vanishing result for the KW equations on a four-manifold $M=M_4$
without boundary.  This will help us know what to aim for.\footnote{The vanishing result without boundary was obtained in \cite{KW}
and the case with a boundary was analyzed in \cite{MW}.}

 The KW equations actually have many different  useful Weitzenbock formulas.  I will first
state some formulas that are useful if we are on a manifold without boundary.   Let $\V=F-\phi\wedge \phi-\star \d_A\phi$, $\W=\d_A\star\phi$,
so the KW equations are $\V=\W=0$.   Clearly then the KW equations are equivalent to the vanishing of 
\be\label{zorf}I=-\int_{M}\Tr\left(\V\wedge\star \V+\W\wedge\star \W\right).\ee
A short calculation gives 
\be\label{morf}I=-\int_{M}\d^4x\sqrt{g}\Tr \left( \frac{1}{2}F_{ij}F^{ij}+D_i\phi_j D^i\phi^j+R_{ij}\phi^i\phi^j+\frac{1}{2}[\phi_i,\phi_j][\phi^i,\phi^j]   \right)\ee
with $R_{ij}$ the Ricci tensor.
If $R_{ij}$ is non-negative, then this is a sum of non-negative terms. The condition $I=0$ forces all these terms to vanish and leads to only a rather trivial class of solutions.

But it is possible to say something useful even if $R_{ij}$ is not non-negative, because
 it is possible to find a family of Weitzenbock formulas.   Define the selfdual 
and anti-selfdual two-forms $\V^+(t)=(F-\phi\wedge\phi+t\d_A\phi)^+$, $\V^-(t)=(F-\phi\wedge\phi-t^{-1}\d_A\phi)^-$.   The equations
$\V^+(t)=\V^-(t)=\W=0$ are a 1-parameter family of elliptic equations, parametrized by $t\in \Bbb{RP}^1$. 
One finds that 
 \begin{align}-\int_M\d^4x \sqrt g&\Tr\left(\frac{t^{-1}}{t+t^{-1}}\V^+_{ij}(t)\V^{+\,ij}(t)+\frac{t}{t+t^{-1}}\V_{ij}^-(t)\V^{-\,ij}(t)+\W^2\right)\cr &
=I +\frac{t-t^{-1}}{4(t+t^{-1})}\int_M \d^4 x \,\epsilon^{ijkl}\Tr F_{ij}F_{kl}. \notag\end{align}

In other words, the same quantity $I$ can be written as a sum of squares in many different ways, modulo the topological invariant
\be\label{noxt}J(t)= \frac{t-t^{-1}}{4(t+t^{-1})}\cdot 32\pi^2 P,~~~ P=\frac{1}{32\pi^2}\int_M \d^4 x \,\epsilon^{ijkl}\Tr F_{ij}F_{kl}.\ee
Now we can deduce the following:  (1) The KW equations cannot have any solutions for $t\not=0,\infty$ except with $P=0$ (if $P\not=0$ for some solution, then by looking at the Weitzenbock formula at some value of $t'$ with $J(t')<J(t)$, we reach a contradiction).  And (2):  If the KW equations are obeyed at one value of $t$ other than $0,\infty$, then they are obeyed at all
$t$.  This is an immediate consequence of the Weitzenbock formula, once we know that $P=0$.   The equations then reduce
to $\mathcal F=0$, where $\mathcal F=\d \A+\A\wedge\A$, with $\A$ the complex connection $\A=A+i\phi$, along with   $\d_A\star\phi=0$. 
According to a theorem of Corlette \cite{Corlette}, the solutions correspond to homomorphisms $\pi_1(M)\to G_\C$ that are in a certain sense semi-stable.

The moral of the story is that the KW equations participate in many different Weitzenbock formulas, not just one, and it is important
to know all of them.   However, none of the formulas that I have written down so far are useful for understanding the Nahm pole
boundary condition.   The reason is that if $\partial M\not=\varnothing$, then the preceding formulas (whose derivation involves integration by parts) will have boundary contributions
if we are on a manifold with boundary, and those boundary contributions
 are divergent in the case of a solution with Nahm pole boundary behavior.    This is inevitable because the expression
that I called $I$ in writing the Weitzenbock formula is divergent in the case of a solution with Nahm pole behavior.  A formula like
\be\label{woxt}-\int\Tr\left(\V\wedge\star\V+\W\wedge \star \W\right)=I+{\mathrm{boundary ~correction}} \ee
must have a boundary correction $-\infty$ in the case of a Nahm pole, since the left hand side is 0 and $I=+\infty$.
A Weitzenbock formula with such divergent terms is not likely to be useful.

To get around this, the best we could hope for would be a Weitenbock formula on $\R^4_+$
in which $I$ is replaced by a sum of squares
of quantities whose vanishing characterizes the Nahm pole solution $A=0$, $\phi=\vec\t\cdot \d\vec x/y$.
 The quantities that vanish\footnote{In the following,  indices $i,j$ take four values corresponding to $x_1,x_2,x_3,y$, but indices $a,b$ take only three values corresponding to $x_1,x_2,x_3$. Also, $\varepsilon_{abc}$ is the Levi-Civita antisymmetric tensor.} in the Nahm pole solution are the curvature $F$, covariant
derivatives and commutators that involve $\phi_y$, namely $D_i\phi_y$ and $[\phi_i,\phi_y]$, covariant derivatives of $\phi_b$
along $\R^3$ such as $D_a\phi_b$, and finally $W_a=D_y\phi_a+\frac{1}{2}\varepsilon_{abc}[\phi_b,\phi_c]$.  (Nahm's equation is $W_a=0$.) 
What we need is true.   Define the following sum of squares of the objects whose vanishing characterizes the Nahm pole solution:
\begin{align}I'=-\int_{\R^3\times \R_+}\d^4x\,\Tr\left(\frac{1}{2}\sum_{i,j}F_{ij}^2+\sum_{a,b}(D_a\phi_b)^2 +\sum_i(D_i\phi_y)^2\right.
\cr\left. +\sum_a[\phi_y,\phi_a]^2
+\sum_a W_a^2\right).\notag\end{align}

Then there is an identity along the lines that we need: 
\be\label{toxt}-\int_M\Tr\left(\V^2+\W^2\right)= I'+\Delta\ee
where $\Delta$ is a certain  boundary term
\be\label{teft}\Delta=-\int_{\partial\R_+^4}\Tr\left(\phi\wedge F +\dots  \right). \ee
(I have omitted some further terms in $\Delta$.)   $\Delta$ is the sum of a contribution at the boundary $y=0$ and on a large hemisphere
$\sqrt{\vec x^2+y^2}>>1$.

Now to get a vanishing theorem that will say that the global Nahm pole solution is the only solution on $\R^4_+$ that
obeys Nahm pole boundary conditions, we need to do the following. We have to prove that if $A,\phi$ approach the Nahm pole
solution for $y\to 0$ and for $\sqrt{\vec x^2+y^2}\to\infty$, then they approach it fast enough so that $\Delta=0$.    Once this is established,  the Weitzenbock formula
will say that a KW solution that satisfies the boundary condition must have $I'=0$.   But $I'$ was constructed so that it vanishes
for and only for the Nahm pole solution.

To find the expected behavior of a solution for $y\to 0$ is a matter of looking at an ODE in which one ignores the $\vec x$
dependence, since that is nonsingular.    In effect, then, we just have to look at the eigenvalues of the linearization of Nahm's
equation (or more exactly a doubled version of Nahm's equation with $A$ as well as $\phi$).    Half of the linearized eigenvalues
are negative and half are positive.  The Nahm pole boundary condition amounts to setting to 0 the coefficients of perturbations with negative
eigenvalues, and allowing the positive ones.  The positive eigenvalues are large enough to ensure that $\Delta=0$ when the Nahm pole
boundary condition is obeyed. 

This shows that there is no contribution to $\Delta$ from the boundary at $y=0$.  To show that there is no contribution to $\Delta$ at $\sqrt{\vec x^2+y^2}\to\infty$, one needs to look at the eigenvalues
of the ``angular'' part of the operator $\L$, which is an operator on a hemisphere $S^3_+$ with Nahm pole boundary conditions along
the boundary.    Those eigenvalues determine how fast a solution will vanish at infinity, assuming that it does vanish at
infinity.  Again the spectrum is such that there is no contribution to $\Delta$. 

This leads to the nonlinear vanishing theorem -- the global Nahm pole solution is the only solution
on the half-space that satisfies the boundary conditions.  Much the same argument proves a linearization of the same statement:
the operator $\L$ obtained by linearizing around the Nahm pole solution has trivial kernel (and hence also trivial cokernel, since $\L$ is
conjugate to $-\L^\dagger$). 

Together with the machinery of uniformly degenerate elliptic operators, this leads to the ellipticity of the Nahm pole boundary
condition in the absence of knots.   But what are we supposed to say in the presence of knots?    As already noted, the knot will be in the boundary  (fig. \ref{thimble}).
To incorporate a knot in the boundary, we  introduce a refinement of the Nahm pole boundary condition, such that $(A,\phi)$ obeys the Nahm pole boundary condition at a generic
boundary point away from a knot, but has some more subtle behavior near the knot.

     \begin{figure}
 \begin{center}
   \includegraphics[width=3in]{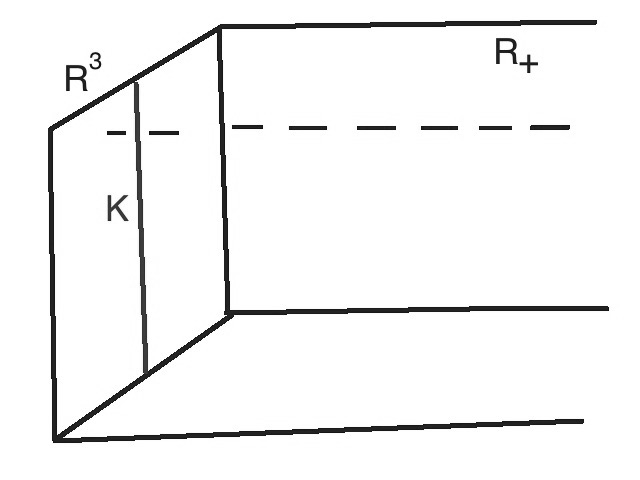}
 \end{center}
\caption{\small A model situation in which a knot is represented by a straight line $K=\R$ in the boundary of a half-space $\R_+^4$.  \label{knotmodel}}
\end{figure}   
  To describe what this more subtle behavior should be, we consider
the case that the knot is locally a straight line $\R\subset \R^3$, so we work on $\R^4_+$ with a knot that lives on a straight line $K$ in the boundary (fig. \ref{knotmodel}).

    \begin{figure}
 \begin{center}
   \includegraphics[width=3in]{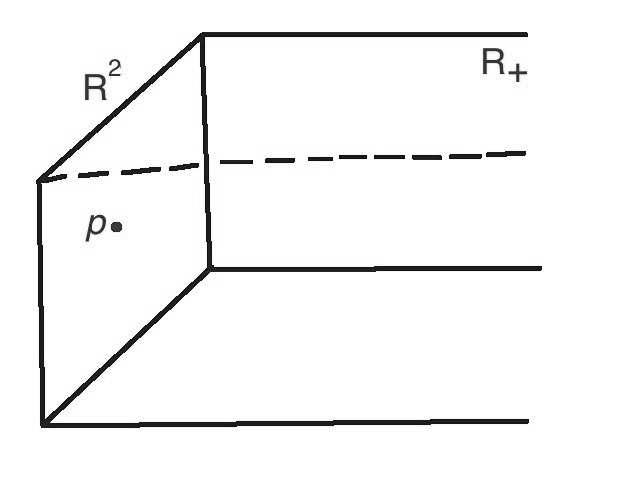}
 \end{center}
\caption{\small Assuming translation invariance along the $K$ direction in fig. \ref{knotmodel}, we reduce to $\R^3_+=\R^2\times \R_+$, with the knot
now represented by a point $p$ in the boundary. \label{Bpoint}}
\end{figure}   

 The idea is going to be to find a singular model solution in the presence of the knot.   This solution will 
 coincide with the Nahm pole
 solution near a boundary point away from $K$, but it will look different near a point of $K$.   The model solution will depend
 on the choice of an irreducible representation $R^\vee$ of the dual group $G^\vee$.    Then a boundary condition is
 defined by saying that one only allows solutions of the KW equations that look like the model solution near a knot.  
 
 For this to make sense, the model solution must look the same near any point of $K$, so we assume that the model
 solution is invariant under translations along $K$.  So we reduce to equations on $\R^2\times \R_+$ with the knot now
 represented by a point $p\in\R^2$ (fig. \ref{Bpoint}).

  Once we reduce to 3 dimensions (and assume vanishing of $A_1$ and $\phi_y$ in a way that can be motivated by the Weitzenbock
  formula) the KW equations become tractable.     Pick coordinates so that $x_1$ runs along the knot $K$; $x_2,x_3$ parametrize
  the normal plane to  $K$ in the boundary; and $y$ measures the distance from the boundary.   
  Define the three operators
\begin{align}\D_1 & = D_2+iD_3=\frac{\partial}{\partial x_2}+i\frac{\partial}{\partial x_3}  +[A_2+iA_3,\,\cdot\, ]\cr
               \D_2& =D_y-i[\phi_1,\,\cdot\,]= \frac{\partial}{\partial y} +[A_y-i\phi_1,\,\cdot\,]\cr 
                  \D_3&=[\phi_2-i\phi_3,\,\cdot\,] , \notag \end{align}  
and also the ``moment map''
\begin{equation} \notag \mu=F_{23}-[\phi_2,\phi_3]-D_y\phi_1.\end{equation}     

The KW equations in this situation become
\be\label{moxy}0=[\D_i,\D_j],~~~i,j=1,\dots,3\ee
along with a ``moment map'' condition
\be\label{loft}\mu=0.\ee
These equations were introduced in \cite{KW} and were called the extended Bogomolny equations.
They are a sort of hybrid of three much-studied equations in the mathematics of gauge theory.    If we drop $\D_1$
(by assuming that the fields are independent of $x_2$ and $x_3$ and that $A_2=A_3=0$), we get Nahm's equation; if we drop $\D_2$
(by assuming that the fields are independent of $y$ and that $A_y=\phi_1=0$) we get Hitchin's equation; and if we drop $\D_3$ (by setting
$\phi_2=\phi_3=0$), we get the Bogomolny equations.

The full system of equations is tractable for the same reason each of those three specializations is.  There are two key facts: 
 (a) the equations $[\D_i,\D_j]=0$
are invariant under $G_\C$-valued gauge transformations ($G_\C$ is the complexification of $G$); 
 (b) the combination of setting $\mu=0$ and dividing by $G$-valued
gauge transformations is equivalent to forgetting the condition $\mu=0$ and dividing by $G_\C$-valued gauge transformations. 
This means that the solutions can be understood in terms of complex geometry.

It is reasonable to expect that the model solution possesses the symmetries of the knot.   So we assume that the model
solution is invariant under a rotation of the boundary $\R^2$ around the point $p\in\R^2$ at which the knot lives, and also invariant
under a scaling of $\R^2\times \R_+$ keeping $p$ fixed.    With these assumptions, the equations $[\D_i,\D_j]=\mu=0$
reduce to affine Toda equations
which are integrable.    One can find all the solutions in closed form, and the solutions that satisfy the Nahm pole boundary condition
away from the knot are classified by an irreducible representation $R^\vee$  of the dual group $G^\vee$.   These solutions were
found for $G=SO(3)$, $G^\vee=SU(2)$ in \cite{Witten} and more generally in \cite{Mikhaylov}. 

How would one go about proving that the KW equations with a boundary condition defined by one of these model equations is
a well-posed (elliptic) problem?    Basically, modulo generalities about uniformly degenerate elliptic operators, we need to
 show that the operator $\L$ obtained by linearizing around one of these
solutions has no kernel or cokernel.   It is again sufficient to show that the kernel vanishes, since $\L^\dagger$ is conjugate to $-\L$.  

Just as in the absence of a knot, we will actually find a nonlinear analog of the vanishing of the kernel of $\L$: any solution of the
KW equations on $\R^3\times \R_+$  (with the knot as an infinite straight line in the boundary, as before) that is asymptotic  to the model solution
both  along the boundary  and at infinity actually coincides
with it. 

The vanishing results we want are the sort that often 
follow from a Weitzenbock formula.    But none of the Weitzenbock formulas that
we considered before are well-adapted to the presence of a knot.   Even the more subtle Weitzenbock formula that includes
the Nahm pole singularity away from a knot
\be\label{tudd}-\int_M\Tr\left(\V^2+\W^2\right)= I'+\Delta\ee does
not give any useful information, because $I'$ (which is the sum of squares of quantities that vanish in the Nahm pole solution without a knot)
is divergent in the presence of a knot so $I'$ will be $+\infty$ and hence $\Delta$ will be $-\infty$ with a knot present.

So we need a new Weitzenbock formula. 
 We  imitate what we did before.    We find a collection of quantities $X_i$ whose vanishing characterizes the model
solution.  (Some obvious $X_i$ are real and imaginary parts of $[\D_i,\D_j]$,  and also $\mu$; the others are quantities like $[\phi_i,\phi_y]$ that vanish because the model
solution has $A_1=\phi_y=0$.)    Then if
\be\label{mudd}I'' =-\sum_i\int_{\R^3\times \R_+}\Tr\,X_i^2,\ee
we have to hope that there is an identity
\be\label{wudd}-\int_{\R^3\times \R_+}\Tr\left(\V\wedge\star \V+\W\wedge\star\W\right)= I''+\widehat \Delta,\ee
where $\widehat \Delta$ is a new boundary term.   It turns out that there is indeed an identity like this.

\begin{figure}
 \begin{center}
   \includegraphics[width=3in]{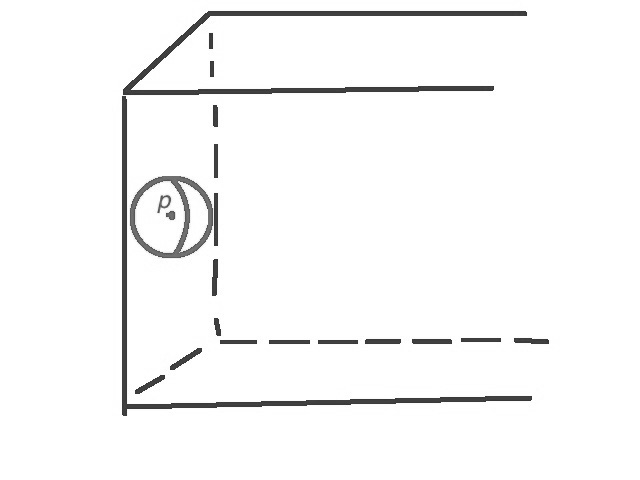}
 \end{center}
\caption{\small To prove the vanishing theorem in the presence of a knot, one has to study the angular operator near a singular point on the boundary. 
The angular operator is defined on the indicated hemisphere. \label{hemi}}
\end{figure}

We still need to show that $\widehat\Delta=0$ 
in the case of a solution that obeys the KW equations and is asymptotic near the knot to the model
solution.   For this, one needs to know what is the asymptotic behavior near the boundary and at infinity of a solution of the KW
equations.     We already know the behavior at a generic boundary point, which was used in our proof of the well-posedness of
the Nahm pole boundary condition without a knot.   To find the behavior near the knot and also at infinity, we now
need to solve the angular part of the equation on a 2d hemisphere $S^2_+$ (fig. \ref{hemi}).

Again it turns out that the eigenvalues of the angular operator are favorable, so there is no contribution to $\widehat \Delta$ either
near $p$ or at infinity.    This fact together with the relevant Weitzenbock formula imply that a solution of the KW equations on
$\R^3\times \R_+$ that is asymptotic on the boundary and at infinity to the model solution with the knot actually coincides with that
model solution.    A linearized version of the same argument shows that the kernel of $\L$ vanishes, which is what we actually
needed to know for ellipticity.

This is the main step in showing that $\L$ is a Fredholm operator in the presence of an arbitrary knot $K$ embedded in
any three-manifold $W$.    Some details are still needed to show that this gives an elliptic boundary condition for the nonlinear
KW equations in the general case of a curved knot \cite{MWtwo}.

\appendix
\section{Some Physics Background}

In this appendix, I will briefly describe some physics background to the treatment of singly-graded Khovanov homology in Lecture One.
Only a bare outline of the string/M-theory context and the  framework of \cite{Witten} for doubly-graded Khovanov homology will be given.
 I aim primarily to explain what is different for the singly-graded theory.

The starting point is the existence of a six-dimensional superconformal field theory with $(0,2)$ supersymmetry, associated to any simply-laced
Lie group $G$, or more precisely to its Dynkin diagram.  This theory has a $\Spin(5)$ group of R-symmetries.   Making use of a subgroup $F=\Spin(2)\subset \Spin(5)$,
the theory can be topologically twisted in such a way that it can be compactified on a Riemann surface $C$ to give a four-dimensional theory with $\mathcal N=2$
supersymmetry.  These are the theories of class S, as studied in  \cite{Gaiotto}.  They have an R-symmetry group $\hat F=(\Spin(3)\times \Spin(2))/\Z_2\cong \mathrm{U}(2)$ (the subgroup of $\Spin(5)$
that commutes with $F$).  Using the $\Spin(3)$ factor, a theory of class S can be topologically twisted and
compactified on a four-manifold $M$ in a way that preserves one supercharge $Q$. This gives
 a theory on $M\times C$ that is topological along
$M$ and holomorphic along $C$.  (The holomorphy along $C$ means, for example, that the cohomology of $Q$ acting on local operators on $M\times C$ varies homomorphically
on $C$.)  The topological-holomorphic theory on $M\times C$ still has an $F=\Spin(2)$ symmetry,  because $F$ commutes with the group $\hat F$ that
was used in the twisting.

The underlying six-dimensional $(0,2)$ model admits half-BPS surface operators that can be supported on any two-manifold $U$ in six dimensions.
However, when the theory is formulated on $M\times C$ as summarized in the last paragraph, if we wish to preserve the supercharge $Q$
of the topological-holomorphic theory, the possible choices of $U$ are quite limited.  $U$ must be of the form $p\times C$, where $p$ is a point in $M$, or  $\Sigma\times q$, where $\Sigma\subset M$ is a two-manifold
and $q$ is a point in $C$.  The reason for this is essentially that any (complete, connected) complex submanifold of $C$ is $C$ itself or a point  $q\in C$.  For
constructing Khovanov homology, we take $U=\Sigma\times q$.

To get singly-graded Khovanov homology, we take $M$ to be simply $\R^4$, and $C$ to be a cylinder $\R\times S^1$.  The six-manifold $M\times C$ is then
simply $\R^4\times \R\times S^1$.  The supercharge $Q$ is invariant under rotations of $\R^4$ (combined with a suitable element of $\Spin(3)\subset \hat F$)
but not under more general rotations of $\R^4\times \R$.  Now we use the fact that the $(0,2)$ model, when formulated on $M'\times S^1$ for any five-manifold $M'$, with
the radius of $M'$ being much greater than that of $S^1$, reduces at long distances on $M'$ 
to maximally supersymmetric gauge theory, with gauge group $G$.  In the context of the topological-holomorphic theory described above, this reduction is valid without taking
any large distance limit.  The resulting supersymmetric gauge theory on $M'=\R^4\times \R$ is infrared-free  (in sharp contrast to the underlying $(0,2)$ model in six dimensions)  and can be
analyzed by classical methods.  In particular, the condition for $Q$-invariance becomes the HW equations, which were mentioned in Lecture One.    These are equations for a pair
consisting of a gauge field $A$, and a field $B$ on $\R^4\times \R$ that is an adjoint-valued section of the pullback to $\R^4\times \R$ of the bundle of selfdual two-forms on $\R^4$.

A surface operator
in six dimensions supported on $\Sigma\times q$ reduces in the gauge theory description to an 't Hooft-like surface operator supported on $\Sigma\times q'$, where $q\in \R\times S^1$
projects to $q'\in \R$.  A solution of the HW equations in the presence of this surface operator is supposed to have a singularity along $\Sigma\times q'$.  
This codimension three singularity should be modeled on the standard codimension three singularity of the Bogomolny equations, suitably embedded in the HW equations with
gauge group $G$.

   To categorify the quantum knot invariants associated to a representation $R^\vee$ of a 
simply-laced\footnote{If $G^\vee$ is not simply-laced, one requires a refinement described in section 5.5 of \cite{Witten}.}   compact Lie
group $G^\vee$, one studies
the HW equations with gauge group $G$ (the Langlands-GNO dual of $G^\vee$, which in particular has the same Lie algebra as $G^\vee$ if $G$ is simply-laced),
and with the appropriate singularity along $\Sigma\times q'$.  The appropriate singularity is obtained, as discussed in Lecture One, by embedding a singular $U(1)$ solution in $G$ using
the homomorphism $\rho:\frak{u}(1)\to\frak t\subset \frak{g}$ that is dual to $R^\vee$.

The HW equations on $\R^4\times \R$ are compatible with the familiar codimension three singularity on a two-manifold  $ V\subset \R^4\times \R$ 
if and only if $V$ is of the form $\Sigma\times q'$ with $\Sigma\subset \R^4$, $q'\in \R$.  This statement is easily verified by inspection of the HW equations.
The explanation that we have given here starting with the $(0,2)$ model serves to explain ``why'' it is true.   The restriction to $V\subset \R^4\times q'\subset \R^4\times \R$
is completely essential for getting knot theory out of this construction, since the relevant topology would disappear if $V$ were free to move in five dimensions.
For example, Khovanov homology arises in the ``time''-independent case $\Sigma=\R_t\times K\subset \R_t\times\R^3=\R^4$, where $K$ is a knot in $\R^3$ and the first factor in $\R^4=\R_t\times \R^3$
is parametrized by the ``time.''  If $V$ were free to vary in $\R^4\times \R$,
then $K$ would be free to vary in $\R^3\times \R$ (the product of the last two factors in $\R^4\times \R=\R_t\times \R^3\times \R$), and could be trivially unknotted.  Much the same point was
made   in Lecture One.   A more general
$V$ (not of the time-independent form $\R_t\times K$) is used to define the ``morphisms'' of Khovanov homology.

What we have described corresponds to the singly-graded version of Khovanov homology (with only the cohomological grading and no ``$q$''-grading).
The single grading comes from the symmetry group $F\cong \Spin(2)$ that was maintained throughout the construction.   To get doubly-graded Khovanov homology,
one takes $C$ to be not $\C^*=\R\times S^1$ but $\C$.  This is done by adding to $\C^*$ a point $q_0$ ``at infinity.''  $\C$ admits an $S^1$ action, leaving fixed the point $q_0$.
In the underlying $(0,2)$ model, one considers a surface operator supported on $U=\Sigma\times q_0$.  Reduction of $M\times C$ on the orbits of $S^1$ leads now to a description
in terms of gauge theory on $M\times \R_+$ (where $\R_+$, a half-line, is the quotient $C/S^1$).  The $S^1$ action leads to the desired second grading.
This doubly-graded version of the construction was the main subject of \cite{Witten}, and the details will not be repeated here.

Research supported in part by NSF Grant PHY-1314311. I thank M. Abouzaid, C. Manolescu, and R. Mazzeo for discussions.

\bibliographystyle{unsrt}

\begin{thebibliography}{99}

\bibitem{GVS}
S. Gukov, A. S. Schwarz, and C. Vafa,
``Khovanov-Rozansky Homology And Topological Strings,''
Lett. Math. Phys. {\bf 74} (2005) 53-74, hep-th/0412243.

\bibitem{OV}
H. Ooguri and C. Vafa, ``Knot Invariants And Topological Strings,''  Nucl. Phys.
{\bf{B577}}
(2000) 419, hep-th/9912123.

\bibitem{Witten}
E. Witten, ``Fivebranes And Knots,'' Quantum Topology {\bf 3} (2012) 1-137, arXiv:1101.3216.

\bibitem{WittenOne}
E. Witten, ``Khovanov Homology and Gauge Theory,''  in  R. Kirby, V. Krushkal, and Z. Wang, eds.,
{\it Proceedings Of The FreedmanFest} (Mathematical Sciences Publishers, 2012) 291-308, arXiv:1108.3103.

\bibitem{WittenTwo}
E. Witten, ``Two Lectures On The Jones Polynomial And Khovanov Homology,'' arXiv:1401.6996.

\bibitem{SS}
P. Seidel and I. Smith, ``A Link Invariant From The Symplectic Geometry of Nilpotent Slices,'' Duke Math. J. {\bf 134}
(2006) 453-514.

\bibitem{Manolescu}
C Manolescu, ``Nilpotent Slices, Hilbert Schemes, and the Jones Polynomial,'' Duke Math. J.  {\bf 132} (2006) 311-369.

\bibitem{AS}
M. Abouzaid and I. Smith, ``Khovanov Homology From Floer Cohomology,'' arXiv:1504.01230.

\bibitem{Jones}
V. F. R. Jones, ``A Polynomial Invariant For Knots Via Von Neumann Algebra,'' Bull. Am. Math. Soc. {\bf 12} (1985) 103-11.

\bibitem{AtiyahFloer}
M. F. Atiyah, ``Floer Homology,'' Progr. Math.  {\bf 133} (1995) 105.

\bibitem{GaiWit}
D. Gaiotto and E. Witten, ``Knot Invariants From Four-Dimensional Gauge Theory,''  Adv. Theor. Math. Phys. {\bf 16} (2012)  935, arXiv:1106.4789.

\bibitem{K} J. Kamnitzer, ``The Beilinson-Drinfeld Grassmannian and Symplectic Knot Homology,'' in D. A. Ellwood (ed.), {\it Grassmannians,
Moduli Spaces, and Vector Bundles} (Clay Math Institute, Cambridge MA and American Mathematical Society, Providence RI, 2011) 81-94.

\bibitem{CK} S. Cautis and J. Kamnitzer, ``Knot Homology Via Derived Categories Of Coherent Sheaves, I. The $\mathfrak{sl}(2)$ Case,''
Duke Math. J. {\bf 142} (2008) 511-588.

\bibitem{thooft} G. 't Hooft, ``On The Phase Transition Towards Permanent Quark Confinement,'' Nucl. Phys. {\bf B138} (1978) 1-25.

\bibitem{KW} A. Kapustin and E. Witten, ``Electric-Magnetic Duality And The Geometric Langlands Program,'' Commun. Numb. Th. Physics
{\bf 1} (2007) 1-236.

\bibitem{Bigelow}
S. Bigelow, ``A Homological Definition Of The Jones Polynomial,'' {\it Geometry and Topology Monographs}, Volume 4: Invariants of Knots and 3-Manifolds
(Kyoto, 2001), pp. 29-41.


\bibitem{Lawrence}
R. Lawrence, ``A Functorial Approach to the One-Variable Jones Polynomial,'' J. Diff. Geom. {\bf 37}  (1993) 689-710.

\bibitem{GalMoore}
D. Galakhov and G. W. Moore, ``Comments On The Two-Dimensional Landau-Ginzburg Approach,''
arXiv:1607.04222.

\bibitem{Haydys}
A. Haydys, ``Fukaya-Seidel Category And Gauge Theory,'' J. Symp. Geom. {\bf 13} (2015) 151-207,  arXiv:1010.2353.

\bibitem{GPV}
S. Gukov, P. Putrov, and C. Vafa, ``Fivebranes and 3-manifold Homology,'' arXiv:1602.05302.


\bibitem{MW}
R. Mazzeo and E. Witten,  ``The Nahm Pole Boundary Condition,''  in {\it The Influence Of Solomon Lefschetz in Geometry and Topology}, ed. L. Katzarkov et. al.
(American Mathematical Society, 2014),
 arXiv:1311.3167.

\bibitem{MWtwo} R. Mazzeo and E. Witten, ``The Nahm Pole Boundary Condition With Knots,'' to appear.



\bibitem{N}
W. Nahm, ``All Selfdual Multimonopoles for Arbitrary Gauge Groups'' (NATO ASI, 1983).

\bibitem{Corlette}
K. Corlette, ``Flat $G$-Bundles With Canonical Metric,'' J. Diff. Geom. {\bf 28} (1988) 361-82.

\bibitem{Mikhaylov}
V. Mikhaylov, ``On The Solutions Of Generalized Bogomolny Equations,'' JHEP 05 (2012) 112, arXiv:1202.4848.

\bibitem{Gaiotto}
D. Gaiotto, ``$\mathcal N=2$ Dualities,'' JHEP {\bf 1208} (2012) 034, arXiv:0904.2715.


\end{thebibliography}

\end{document}